\documentclass[a4paper,11pt,twoside,english]{article}
\usepackage{babel}
\usepackage{latexsym,amsfonts}
\usepackage{amsthm}
\usepackage{amsmath}
\usepackage{mathtools}
\usepackage{paralist}

\usepackage{lastpage}
\usepackage{fancyhdr}
\fancypagestyle{plain}{

  \fancyhead[L,C]{}
  \fancyhead[R]{\today}
  \fancyfoot[R,L]{}
  \fancyfoot[C]{\thepage \ of \pageref{LastPage}}
}

\title{On certain geometric properties in Banach spaces of vector-valued functions}
\author{Jan-David Hardtke}
\date{}

\DeclareMathOperator{\aco}{aco}

\providecommand{\sm}{\setminus}
\providecommand{\ssq}{\subseteq}

\providecommand{\N}{\ensuremath{\mathbb{N}}}

\providecommand{\R}{\ensuremath{\mathbb{R}}}

\providecommand{\K}{\ensuremath{\mathbb{K}}}
\providecommand{\T}{\ensuremath{\mathbb{T}}}
\providecommand{\E}{\ensuremath{\mathcal{E}}}
\providecommand{\A}{\ensuremath{\mathcal{A}}}
\providecommand{\U}{\ensuremath{\mathcal{U}}}
\providecommand{\eps}{\ensuremath{\varepsilon}}
\providecommand{\dist}[2]{\ensuremath{\operatorname{d} (#1,#2)}}

\providecommand{\keywords}[1]{
{\let\thefootnote=\relax
\footnote{{\em Keywords}: #1}}
\addtocounter{footnote}{-1}
}

\providecommand{\AMS}[1]{
{\let\thefootnote=\relax
\footnote{{\em AMS Subject Classification} (2010): #1}}
\addtocounter{footnote}{-1}
}

\providecommand{\address}{
{\sc \noindent Department of Mathematics \\
Freie Universit\"at Berlin \\
Arnimallee 6, 14195 Berlin \\
Germany \\}
}

\DeclarePairedDelimiter{\set}{\lbrace}{\rbrace}
\DeclarePairedDelimiter{\paren}{\lparen}{\rparen}

\DeclarePairedDelimiter{\abs}{\lvert}{\rvert}
\DeclarePairedDelimiter{\norm}{\lVert}{\rVert}

\theoremstyle{definition}
\newtheorem{definition}{Definition}[section]
\newtheorem*{definition*}{Definition}

\newtheorem*{example*}{Example}

\newtheorem*{remark*}{Remark}
\theoremstyle{plain}
\newtheorem{lemma}[definition]{Lemma}
\newtheorem*{lemma*}{Lemma}
\newtheorem{proposition}[definition]{Proposition}
\newtheorem*{proposition*}{Proposition}
\newtheorem{theorem}[definition]{Theorem}
\newtheorem*{theorem*}{Theorem}
\newtheorem{corollary}[definition]{Corollary}
\newtheorem*{corolary*}{Corollary}

\newenvironment{Proof}[1][\proofname]{\begin{proof}[#1] \setlength{\parindent}{0pt}}{\end{proof}}
\newenvironment{Abstract}{\centering\begin{minipage}{0.8\textwidth} \noindent \small {\sc Abstract.}}{\end{minipage}\par}

\usepackage{color}
\definecolor{darkgreen}{rgb}{0,0.5,0}

\numberwithin{equation}{section}
\newtagform{colored}[\color{blue}]{\color{blue}(}{\color{blue})}
\usetagform{colored}

\usepackage[colorlinks,linkcolor=blue,citecolor=red,urlcolor=darkgreen]{hyperref}
\providecommand{\email}{{\it E-mail address:} \href{mailto:hardtke@math.fu-berlin.de}{\tt hardtke@math.fu-berlin.de}}

\usepackage{amsrefs}

\begin{document}

\maketitle

\begin{Abstract}
We consider a certain type of geometric properties of Banach spaces, which includes
for instance octahedrality, almost squareness, lushness and the Daugavet property. For this type of properties, 
we obtain a general reduction theorem, which, roughly speaking, states the following: if the 
property in question is stable under certain finite absolute sums (for example finite $\ell^p$-sums), 
then it is also stable under the formation of corresponding K\"othe-Bochner spaces (for example $L^p$-Bochner spaces).\par
From this general theorem, we obtain as corollaries a number of new results as well as some alternative proofs of
already known results concerning octahedral and almost square spaces and their relatives, diameter-two-properties, lush spaces 
and other classes.
\end{Abstract}
\AMS{46B20 46E40}

\section{Introduction}\label{sec:intro}
Let $\K$ be the real or complex field. We consider a class $\E$ of Banach spaces over $\K$
which is closed under isometric isomorphisms, i.\,e. if $X\in \E$ and $Y$ is isometrically
isomorphic to $X$, then $Y\in \E$.\par 
For a given Banach space $X$, we denote by $X^*$ its dual space, by $B_X$ its closed unit ball 
and by $S_X$ its unit sphere. Furthermore, $B_X^{\text{fin}}$ and $S_X^{\text{fin}}$ will denote 
the sets of all finite sequences in $B_X$ and $S_X$. For fixed $n\in \N$, $B_X^n$ and $S_X^n$ will
stand for the sets of all sequences of length $n$ in $B_X$ and $S_X$. Given $\mathbf{x}=(x_1,\dots,x_n)\in B_X^n$, 
we set $\norm{\mathbf{x}}_{\infty}:=\max_{i=1,\dots,n}\norm{x_i}$. Finallly, $\U(X)$ denotes the set of all closed, nontrivial subspaces of $X$.\par
The following is our main definition.
\begin{definition}\label{def:test}
Let $X$ be a Banach space. A family of real-valued functions 
$F_{\eps,U}$ on $B_U^{\text{fin}}\times B_{U^*}^{\text{fin}}\times B_U^{\text{fin}}\times B_{U^*}^{\text{fin}}$
with $U\in \U(X)$ and $\eps>0$ is said to be a test family for $\E$ in $X$ if the following conditions are satisfied:
\begin{enumerate}[\upshape(i)]
\item For every $U\in \U(X)$ one has that $U\in \E$ if and only if for every $\eps>0$ and all $\mathbf{x}\in S_U^{\text{fin}}$
and $\mathbf{x}^*\in S_{U^*}^{\text{fin}}$ there exist $\mathbf{y}\in S_U^{\text{fin}}$ and $\mathbf{y}^*\in S_{U^*}^{\text{fin}}$ 
such that $F_{\eps,U}(\mathbf{x},\mathbf{x^*},\mathbf{y},\mathbf{y}^*)\leq \eps$.
\item If $0<\eps_1<\eps_2$ and $U\in \U(X)$, then $F_{\eps_1,U}\geq F_{\eps_2,U}$.			
\item There exists $c>0$ such that for all $U\in \U(X)$, all $\eps>0$, every $\mathbf{x},\mathbf{y}\in B_U^{\text{fin}}$ and  
every $\mathbf{x}^*, \mathbf{y}^*\in B_{X^*}^{\text{fin}}$ one has
\begin{equation*}
F_{\eps,X}(\mathbf{x},\mathbf{x}^*,\mathbf{y},\mathbf{y}^*)\leq cF_{\eps,U}(\mathbf{x},\mathbf{x}^*|_U,\mathbf{y},\mathbf{y}^*|_U),			
\end{equation*}
where for $\mathbf{x}^*=(x_1^*,\dots,x_n^*)$ we define $\mathbf{x}^*|_U=(x_1^*|_U,\dots,x_n^*|_U)$ (and analogously for $\mathbf{y}^*$).
\item For every $\eps>0$, all $\tau>0$, each $\mathbf{x}^*\in B_{X^*}^{\text{fin}}$, all $n\in \N$ and all $\mathbf{x}\in B_X^n$ there exists 
a $\delta>0$ such that 
\begin{equation*}
|F_{\varepsilon,X}(\mathbf{x},\mathbf{x}^*,\mathbf{y},\mathbf{y}^*)-F_{\varepsilon,X}(\mathbf{z},\mathbf{x}^*,\mathbf{y},\mathbf{y}^*)|\leq \tau
\end {equation*}
holds for all $\mathbf{y}\in B_X^{\text{fin}}$, all $\mathbf{y}^*\in B_{X^*}^{\text{fin}}$ and every $\mathbf{z}\in B_X^n$ with 
$\norm{\mathbf{x}-\mathbf{z}}_{\infty}\leq \delta$.
\item For every $\eps>0$, all $n,m\in \N$ and all $\eta>0$ there exists $\theta>0$ such that for every $U\in \U(X)$ one has 
\begin{equation*}
|F_{\varepsilon,U}(\mathbf{x},\mathbf{x}^*,\mathbf{y},\mathbf{y}^*)-F_{\varepsilon,U}(\mathbf{x},\mathbf{z}^*,\mathbf{y},\mathbf{y}^*)|\leq \eta
\end{equation*}
for all $\mathbf{x}\in B_U^n$, all $\mathbf{y}\in B_U^{\text{fin}}$, all $\mathbf{y}^*\in B_{U^*}^{\text{fin}}$ and all $\mathbf{x}^*, \mathbf{z}^*\in B_{U^*}^m$
with $\norm{\mathbf{x}^*-\mathbf{z}^*}_{\infty}\leq \theta$.
\end{enumerate}
\end{definition}

Roughly speaking, we want to show that if a Banach space property can be characterised in terms of test families
and is stable under certain finite, absolute sums, then it is also stable under the formation of corresponding 
K\"othe-Bochner function spaces.\par 
Examples of Banach space properties which can be described by test families will be presented in the next section
(the constant $c$ in the above definition will be 1 for all these examples). Here we continue with the necessary 
basics on absolute sums and K\"othe-Bochner spaces.\par 
Let $I$ be a non-empty set, $E$ a subspace of $\R^I$ with $e_i\in E$ for all $i\in I$ and $\norm*{\,.\,}_E$ a 
complete norm on $E$ (here $e_i$ denotes the characteristic function of $\set*{i}$).\par
\noindent The norm $\norm*{\,.\,}_E$ is called {\em absolute} if 
\begin{align*}
&(a_i)_{i\in I}\in E, \ (b_i)_{i\in I}\in \R^I \ \mathrm{and} \ \abs*{a_i}=\abs*{b_i} \ \forall i\in I \\
&\Rightarrow \ (b_i)_{i\in I}\in E \ \mathrm{and} \ \norm*{(a_i)_{i\in I}}_E=\norm*{(b_i)_{i\in I}}_E.
\end{align*}
The norm is called {\em normalised} if $\norm*{e_i}_E=1$ for every $i\in I$.\par
Standard examples of subspaces of $\R^I$ with absolute normalised norm are of course the spaces $\ell^p(I)$
for $1\leq p\leq \infty$ and the space $c_0(I)$.\par
We note the following lemma on absolute norms (see e.\,g. \cite{lee}*{Remark 2.1}).
\begin{lemma}\label{lemma:abs norms}
Let $(E,\norm*{\,.\,}_E)$ be a subspace of $\R^I$ with an absolute normalised norm. Then the following is true.
\begin{align*}
&(a_i)_{i\in I}\in E, \ (b_i)_{i\in I}\in \R^I \ \mathrm{and} \ \abs*{b_i}\leq\abs*{a_i} \ \forall i\in I \\
&\Rightarrow \ (b_i)_{i\in I}\in E \ \mathrm{and} \ \norm*{(b_i)_{i\in I}}_E\leq\norm*{(a_i)_{i\in I}}_E.
\end{align*}
\end{lemma}
If $(X_i)_{i\in I}$ is a family of (real or complex) Banach spaces we put
\begin{equation*}
\Bigl[\bigoplus_{i\in I}X_i\Bigr]_E:=\set*{(x_i)_{i\in I}\in \prod_{i\in I}X_i: (\norm*{x_i})_{i\in I}\in E}.
\end{equation*}
This defines a subspace of the product space $\prod_{i\in I}X_i$ which becomes a Banach space when endowed with the norm 
\begin{equation*}
\norm*{(x_i)_{i\in I}}_E:=\norm*{(\norm*{x_i})_{i\in I}}_E \ \forall (x_i)_{i\in I}\in \Bigl[\bigoplus_{i\in I}X_i\Bigr]_E.
\end{equation*}
We call this Banach space the absolute sum of the family $(X_i)_{i\in I}$ with respect to $E$. For $p\in [1,\infty]$ and 
$E=\ell^p(I)$ one obtains the usual $p$-sums of Banach spaces.\par 
The ``continuous counterpart'' to absolute sums are the K\"othe-Bochner function spaces, whose definition we will recall now. 
Let $(S,\A,\mu)$ be a complete, $\sigma$-finite measure space. For $A\in \A$ we denote by $\chi_A$ the characteristic function 
of $A$. A K\"othe function space over $(S,\A,\mu)$ is a Banach space $(E,\norm{\cdot}_E)$ of real-valued measurable 
functions on $S$ (modulo equality $\mu$-almost everywhere) such that 
\begin{enumerate}[(i)]
\item $\chi_A\in E$ for every $A\in \A$ with $\mu(A)<\infty$,
\item for every $f\in E$ and every set $A\in \A$ with $\mu(A)<\infty$ $f$ is $\mu$-integrable over $A$,
\item if $g$ is measurable and $f\in E$ such that $\abs*{g(t)}\leq\abs*{f(t)}$ $\mu$-a.\,e. then $g\in E$
and $\norm{g}_E\leq\norm{f}_E$.
\end{enumerate}
Standard examples are the spaces $L^p(\mu)$ for $1\leq p\leq\infty$.\par
Further recall that, given a Banach space $X$, a function $f:S \rightarrow X$ is called simple if there are finitely many disjoint 
measurable sets $A_1,\dots ,A_n\in \A$ such that $\mu(A_i)<\infty$ for all $i=1,\dots,n$, $f$ is constant on each $A_i$ and $f(t)=0$
for every $t\in S\sm \bigcup_{i=1}^nA_i$. The function $f$ is said to be Bochner-measurable if there exists a sequence $(f_n)_{n\in \N}$ 
of simple functions such that $\lim_{n\to \infty}\norm{f_n(t)-f(t)}=0$ $\mu$-a.\,e.\par 
For a K\"othe function space $E$ and a Banach space $X$, we denote by $E(X)$ the space of all Bochner-measurable functions
$f:S\rightarrow X$ (modulo equality a.\,e.) such that $\norm{f(\cdot)}\in E$. Endowed with the norm $\norm{f}_{E(X)}=
\norm*{\norm{f(\cdot)}}_E$ $E(X)$ becomes a Banach space, the so called K\"othe-Bochner space induced by $E$ and $X$.
For $E=L^p(\mu)$ we obtain the usual Lebesgue-Bochner spaces $L^p(\mu,X)$ for $1\leq p\leq\infty$. For more information on 
K\"othe-Bochner spaces the reader is referred to the book \cite{lin}.

\section{Examples}\label{sec:examples}
We will now discuss a number of examples of Banach space properties which can be described via test families.
We start with the octahedral spaces and their relatives.

\subsection{Octahedrality}\label{sub:octa}
A real Banach space $X$ is called octahedral (OH) (see \cite{godefroy}) if the following holds: for every 
finite-dimensional subspace $F$ of $X$ and every $\eps>0$ there is some $y\in S_X$ such that
\begin{equation*}
\norm{x+y}\geq (1-\eps)(\norm{x}+1) \ \ \forall x\in F.
\end{equation*}
$\ell^1$ is the model example of an OH space. It is known that a Banach space has an equivalent OH norm if and only if 
it contains an isomorphic copy of $\ell^1$ (see \cite{deville}*{Theorem 2.5, p. 106}).\par 
In the paper \cite{haller3}, two variants of octahedrality where introduced.\par 
$X$ is called locally octahedral (LOH) if for every $x\in X$ and every $\eps>0$ there exists $y\in S_X$ such that
\begin{equation*}
\norm{sx+y}\geq (1-\eps)(|s|\norm{x}+1) \ \ \forall s\in \R.
\end{equation*}
$X$ is called weakly octahedral (WOH) if for every finite-dimensional subspace $F$ of $X$, every $x^*\in B_{X^*}$ 
and each $\eps>0$ there is some $y\in S_X$ such that 
\begin{equation*}
\norm{x+y}\geq (1-\eps)(|x^*(x)|+1) \ \ \forall x\in F.
\end{equation*}
The motivation for this definition in \cite{haller3} was the study of so called diameter-two-properties.
Given $x^*\in S_{X^*}$ and $\alpha>0$, the slice of $B_X$ induced by $x^*$ and $\alpha$ is
$S(x^*,\alpha):=\set*{z\in B_X,x^*(z)>1-\alpha}$. According to \cite{abrahamsen}, the space $X$ is said 
to have the local diameter-two-property (LD2P) if every slice of $B_X$ has diamter 2; 
$X$ has the diameter-two-property (D2P) if every nonempty, relatively weakly open subset of $B_X$ has diameter 2; 
$X$ has the strong diameter-two-property (SD2P) if every convex combination of slices of $B_X$ has diameter 2.\par 
The following results were proved in \cite{haller3}: 
\begin{enumerate}[(a)]
\item $X$ has LD2P $\iff$ $X^*$ is LOH.
\item $X$ has D2P $\iff$ $X^*$ is WOH.
\item $X$ has SD2P $\iff$ $X^*$ is OH.
\end{enumerate}
The equivalence (c) was also proved independently in \cite{becerra-guerrero2}.\par
It is known that the three diameter-two-properties are really different. For example, it follows from the 
results on direct sums in \cite{haller3} that $c_0\oplus_2 c_0$ has the D2P but not the SD2P (we will
recall these results in Section \ref{sec:appl}).\par
Concerning the nonequivalence of the LD2P and the D2P, it has been shown in \cite{becerra-guerrero4} that there is 
a Banach space with the LD2P whose unit ball contains relatively weakly open subsets of arbitrarily small diameter
(every Banach space containing an isomorphic copy of $c_0$ can be renormed to become such a space 
\cite{becerra-guerrero4}*{Theorem 2.4}).\footnote{Note that the abbreviation SD2P in \cite{becerra-guerrero4} does 
not stand for ``strong diameter-two-property'' but for ``slice diameter-two-property'', which coincides 
with the LD2P of \cite{abrahamsen}.}\par 
In \cite{kubiak} it was shown that Ces\`aro function spaces have the D2P.\par
It is possible to characterise all three octahedrality properties in terms of test families. To do so,
we make use of the following equivalent formulations proved in \cite{haller3} (other equivalent 
characterisations in terms of coverings of the unit ball were proved in \cite{haller2}).\par 
A Banach space $X$ is OH if and only if for every $n\in \N$, all 
$x_1,\dots,x_n\in S_X$ and every $\eps>0$ there exists an element $y\in S_X$
such that $\norm{x_i+y}\geq 2-\eps$ for all $i=1,\dots,n$.\par
$X$ is LOH if and only if for every $x\in S_X$ and all $\eps>0$
there exists $y\in S_X$ such $\norm{x\pm y}\geq 2-\eps$.\par
Of course, the same characterisations also hold for all closed subspaces of $X$.
Thus if we put
\begin{equation*} 
F_{\eps,U}(\mathbf{x},\mathbf{x}^*,\mathbf{y},\mathbf{y}^*):=
\max\set*{2-\norm{x_i+y_1}:i=1,\dots,n}
\end{equation*}
for $U\in \U(X)$, $\mathbf{x}=(x_1,\dots,x_n),\mathbf{y}=(y_1,\dots,y_m)\in B_U^{\text{fin}}$ 
and $\mathbf{x}^*,\mathbf{y}^*\in B_{U^*}^{\text{fin}}$,
we obtain a test family for the class of octahedral spaces in $X$.\par
If we put instead 
\begin{equation*} 
F_{\eps,U}(\mathbf{x},\mathbf{x}^*,\mathbf{y},\mathbf{y}^*):=
\max\set*{2-\norm{x_1+y_1},2-\norm{x_1-y_1}},
\end{equation*}
we obtain a test family for LOH in $X$.\par
(In both cases, condition (i) in Definition \ref{def:test} follows from the above characterisations, 
while the conditions (ii)--(v) are easily verified.)\par 
For weak octahedrality, the following was proved in \cite{haller3}: $X$ is WOH if and only if 
for every $n\in \N$, all $x_1,\dots,x_n\in S_X$ , every $x^*\in S_{X^*}$ and every $\eps>0$ there 
exists a $y\in S_X$ such that $\norm{x_i+ty}\geq (1-\eps)(|x^*(x_i)|+t)$ for all $i=1,\dots,n$ and 
every $t\geq \eps$.\par
(The original formulation in \cite{haller3} reads ``for every $x^*\in B_{X^*}$'', but it clearly 
suffices to take $x^*\in S_{X^*}$.)\par 
Thus if we define 
\begin{equation*} 
F_{\eps,U}(\mathbf{x},\mathbf{x}^*,\mathbf{y},\mathbf{y}^*):=
\max_{i=1,\dots,n}\sup_{t\geq \eps}\paren*{1-\frac{\norm{x_i+ty_1}}{|x_1^*(x_i)|+t}}
\end{equation*}
for $U\in \U(X)$, $\mathbf{x}=(x_1,\dots,x_n),\mathbf{y}=(y_1,\dots,y_m)\in B_U^{\text{fin}}$ 
and $\mathbf{x}^*=(x_1^*,\dots,x_k^*)$, $\mathbf{y}^*\in B_{U^*}^{\text{fin}}$, then 
condition (i) in the definition of a test family for the class of WOH spaces is satisfied and (ii)
and (iii) are clearly true as well. The conditions (iv) and (v) easily follow from the next auxiliary lemma.\par
\begin{lemma}\label{lemma:aux}
If $Y$ is a real Banach space and $\eps>0$, define the function $f:B_Y\times B_Y\times B_{Y^*} \rightarrow \R$ by
\begin{equation*}
f(x,y,x^*):=\sup_{t\geq \eps}\paren*{1-\frac{\norm{x+ty}}{|x^*(x)|+t}} \ \ \forall x,y\in B_Y, \forall x^*\in B_{Y^*}.
\end{equation*}
If $\delta>0$, $x,\tilde{x},y,\tilde{y}\in B_Y$ with $\norm{x-\tilde{x}}, \norm{y-\tilde{y}}\leq \delta$ and
$x^*,\tilde{x}^*\in B_{Y^*}$ with $\norm{x^*-\tilde{x}^*}\leq \delta$, then
\begin{equation*}
|f(x,y,x^*)-f(\tilde{x},\tilde{y},\tilde{x}^*)|\leq \delta (3/\eps+2/\eps^2+1).
\end{equation*}
\end{lemma}

\begin{Proof}
We have $||\tilde{x}^*(\tilde{x})|-|x^*(x)||\leq |\tilde{x}^*(\tilde{x})-\tilde{x}^*(x)|+|\tilde{x}^*(x)-x^*(x)|\leq 2\delta$.\par 
Thus, for every $t\geq \eps$ we have
\begin{align*}
&1-\frac{\norm{\tilde{x}+t\tilde{y}}}{|\tilde{x}^*(\tilde{x})|+t}-\paren*{1-\frac{\norm{x+ty}}{|x^*(x)|+t}} \\
&=\frac{\norm{x+ty}(|\tilde{x}^*(\tilde{x})|+t)-\norm{\tilde{x}+t\tilde{y}}(|x^*(x)|+t)}{(|x^*(x)|+t)(|\tilde{x}^*(\tilde{x})|+t)} \\
&\leq \frac{\norm{x+ty}(|x^*(x)|+2\delta+t)-\norm{\tilde{x}+t\tilde{y}}(|x^*(x)|+t)}{(|x^*(x)|+t)(|\tilde{x}^*(\tilde{x})|+t)} \\
&\leq \frac{2\delta}{t^2}\norm{x+ty}+\frac{\norm{x+ty}-\norm{\tilde{x}+t\tilde{y}}}{|\tilde{x}^*(\tilde{x})|+t} \\
&\leq \frac{2\delta}{t^2}(1+t)+\frac{\norm{x-\tilde{x}}+t\norm{y-\tilde{y}}}{|\tilde{x}^*(\tilde{x})|+t} \\
&\leq \frac{2\delta}{t^2}(1+t)+\frac{(1+t)\delta}{t}=\delta(3/t+2/t^2+1) \\
&\leq \delta (3/\eps+2/\eps^2+1).
\end{align*}
By symmetry we also have
\begin{equation*}
1-\frac{\norm{x+ty}}{|x^*(x)|+t}-\paren*{1-\frac{\norm{\tilde{x}+t\tilde{y}}}{|\tilde{x}^*(\tilde{x})|+t}}\leq \delta (3/\eps+2/\eps^2+1)
\end{equation*}
for all $t\geq \eps$. This implies the desired inequality.
\end{Proof}

The above-mentioned dual characterisations from \cite{haller3} allow us to write the diameter-two-properties 
in terms of test families as well.\par 
Since a Banach space has the LD2P if and only if its dual is LOH, we obtain a test family for the LD2P in $X$ by setting
\begin{equation*} 
F_{\eps,U}(\mathbf{x},\mathbf{x}^*,\mathbf{y},\mathbf{y}^*):=
\max\set*{2-\norm{x_1^*+y_1^*},2-\norm{x_1^*-y_1^*}}
\end{equation*}
for $U\in \U(X)$, $\mathbf{x},\mathbf{y}\in B_U^{\text{fin}}$ and 
$\mathbf{x}^*=(x_1^*,\dots,x_n^*)$, $\mathbf{y}^*=(y_1^*,\dots,y_m^*)\in B_{U^*}^{\text{fin}}$.\par
Likewise, since a Banach has the SD2P if and only its dual is OH, a test family for the SD2P in $X$ is given by
\begin{equation*} 
F_{\eps,U}(\mathbf{x},\mathbf{x}^*,\mathbf{y},\mathbf{y}^*):=
\max\set*{2-\norm{x_i^*+y_1^*}:i=1,\dots,n}
\end{equation*}
for $U\in \U(X)$, $\mathbf{x},\mathbf{y}\in B_U^{\text{fin}}$ and 
$\mathbf{x}^*=(x_1^*,\dots,x_n^*)$, $\mathbf{y}^*=(y_1^*,\dots,y_m^*)\in B_{U^*}^{\text{fin}}$.
(In both cases the conditions (i)--(v) are easily checked.)\par 
We also know that a Banach space has the D2P if and only if its dual is WOH. Then we can make use of 
the following characterisation (see \cite{haller3}) for the property WOH in dual spaces, which does 
not involve the bidual:\par
$X^*$ is WOH if and only if for every $n\in \N$, all $x_1^*,\dots,x_n^*\in S_{X^*}$, every $x\in S_X$ 
and every $\eps>0$ there exists $y^*\in S_{X^*}$ such that 
\begin{equation*}
\norm{x_i^*+ty^*}\geq (1-\eps)(|x_i^*(x)|+t) \ \ \forall i\in \set*{1,\dots,n}\ \forall t\geq \eps.
\end{equation*}
Thus we can define a test family for the D2P in $X$ by 
\begin{equation*} 
F_{\eps,U}(\mathbf{x},\mathbf{x}^*,\mathbf{y},\mathbf{y}^*):=
\max_{i=1,\dots,n}\sup_{t\geq \eps}\paren*{1-\frac{\norm{x_i^*+ty_1^*}}{|x_i^*(x_1)|+t}}
\end{equation*}
for $U\in \U(X)$, $\mathbf{x}=(x_1,\dots,x_k),\mathbf{y}\in B_U^{\text{fin}}$ 
and $\mathbf{x}^*=(x_1^*,\dots,x_n^*)$, $\mathbf{y}^*=(y_1^*,\dots,y_m^*)\in B_{U^*}^{\text{fin}}$.\par 
(The conditions (i)--(iii) are clear, and the conditions (iv) and (v) are proved by using an auxiliary lemma 
similar to Lemma \ref{lemma:aux}).\par 
We remark that it is also possible to describe the LD2P via a different test family, using directly the definition of the LD2P
(and not its dual characterisation). It is easily checked that a Banach space $X$ has the LD2P if and only if the following holds: 
for every $x^*\in S_{X^*}$ and every $\eps>0$ there exist $y_1,y_2\in S_X$ such that $x^*(y_1),x^*(y_2)\geq 1-\eps$ and 
$\norm{y_1-y_2}\geq 2-\eps$.\par 
Thus we can define a test family for the LD2P in $X$ as follows:
\begin{equation*} 
F_{\eps,U}(\mathbf{x},\mathbf{x}^*,\mathbf{y},\mathbf{y}^*):=
\max\set*{1-x_1^*(y_1),1-x_1^*(y_2),2-\norm{y_1-y_2}}
\end{equation*}
for $U\in \U(X)$, $\mathbf{x},\mathbf{y}=(y_1,\dots,y_m)\in B_U^{\text{fin}}$ 
and $\mathbf{x}^*=(x_1^*,\dots,x_n^*),\mathbf{y}^*\in B_{U^*}^{\text{fin}}$, where $y_2:=y_1$
if $m=1$ (once again, the conditions (i)--(v) in Definition \ref{def:test} are easily verified).\par 
Finally, there is yet another weakening of the definition of octahedral spaces, which was introduced in \cite{haller4}:
$X$ is called alternatively octahedral (AOH) if for every $n\in \N$, all $x_1,\dots,x_n\in S_X$ and 
every $\eps>0$ there is some $y\in S_X$ such that 
\begin{equation*}
\max\set*{\norm{x_i+y},\norm{x_i-y}}\geq 2-\eps \ \ \forall i=1,\dots,n.
\end{equation*}
Every octahedral space is alternatively octahedral, while for example $c_0$ is alternatively 
octahedral but not locally octahedral (see \cite{haller4}).\par
It is easily checked that 
\begin{equation*} 
F_{\eps,U}(\mathbf{x},\mathbf{x}^*,\mathbf{y},\mathbf{y}^*):=
\max\set*{2-\max\set*{\norm{x_i+y_1},\norm{x_i-y_1}}:i=1,\dots,n},
\end{equation*}
where $U\in \U(X)$, $\mathbf{x}=(x_1,\dots,x_n),\mathbf{y}=(y_1,\dots,y_m)\in B_U^{\text{fin}}$ 
and $\mathbf{x}^*,\mathbf{y}^*\in B_{U^*}^{\text{fin}}$, defines a test family for AOH in $X$.

\subsection{Almost square spaces}\label{sub:asq}
Next we turn to the classes of almost square and locally almost square Banach spaces. These notions 
were introduced in \cite{abrahamsen2}.\par 
A real Banach space $X$ is said to be almost square (ASQ) if the following holds: for all $n\in \N$ 
and all $x_1,\dots,x_n\in S_X$ there exists a sequence $(y_k)_{k\in \N}$ in $B_X$ such that
$\norm{y_k}\to 1$ and $\norm{x_i+y_k}\to 1$ for all $i=1,\dots,n$.\par 
$X$ is called locally almost square (LASQ) if for every $x\in S_X$ there is a sequence $(y_k)_{k\in \N}$ 
in $B_X$ such that $\norm{y_k}\to 1$ and $\norm{x\pm y_k}\to 1$.\par 
According to \cite{abrahamsen2} $X$ is ASQ if and only if for every $\eps>0$, every $n\in \N$ and all
$x_1,\dots,x_n\in S_X$ there exists a $y\in S_X$ such that $\norm{x_i-y}\leq 1+\eps$ for all $i=1,\dots,n$,
and $X$ is LASQ if and only if for every $\eps>0$ and every $x\in S_X$ there is some $y\in S_X$ such that
$\norm{x\pm y}\leq 1+\eps$.\par 
$c_0$ is the model example of an ASQ space. It was further proved in \cite{abrahamsen2} that every ASQ space 
contains an isomorphic copy of $c_0$ and that every separable Banach spaces containing an isomorphic copy of 
$c_0$ has an equivalent ASQ norm. In \cite{becerra-guerrero3} it was proved that the same holds also for 
nonseparable spaces.\par 
In \cite{abrahamsen2} it was also proved that $X^*$ is OH (i.\,e. $X$ has the SD2P) whenever $X$ is ASQ.
By \cite{kubiak}*{Proposition 2.5} every LASQ space has the LD2P.\par 
If we define 
\begin{equation*} 
F_{\eps,U}(\mathbf{x},\mathbf{x}^*,\mathbf{y},\mathbf{y}^*):=
\max\set*{\norm{x_i-y_1}-1:i=1,\dots,n}
\end{equation*}
for $U\in \U(X)$, $\mathbf{x}=(x_1,\dots,x_n),\mathbf{y}=(y_1,\dots,y_m)\in B_U^{\text{fin}}$ and 
$\mathbf{x}^*,\mathbf{y}^*\in B_{U^*}^{\text{fin}}$, then we obtain a test family for ASQ in $X$, 
as is easily checked.\par 
Likewise, a test family for LASQ in $X$ is given by 
\begin{equation*} 
F_{\eps,U}(\mathbf{x},\mathbf{x}^*,\mathbf{y},\mathbf{y}^*):=
\max\set*{\norm{x_1+y_1}-1,\norm{x_1-y_1}-1}.
\end{equation*}
There is also an intermediate notion of weakly almost square (WASQ) spaces defined in \cite{abrahamsen2}
(by \cite{kubiak}*{Proposition 2.6} theses spaces have the D2P) but it is not clear whether this notion 
can be phrased in terms of test families.

\subsection{The Daugavet property}\label{sub:daugavet}
We now consider spaces with the Daugavet and the alternative Daugavet property.\par 
A real Banach space $X$ is said to have the Daugavet property (DP) if the equality 
$\norm{\text{id}+T}=1+\norm{T}$ holds for every rank-one operator $T:X \rightarrow X$
(see for example \cites{kadets, werner}).\par
Examples of such spaces include $C(K)$ for compact Hausdorff spaces $K$ without isolated points, and
$L^1(\mu)$ for atomless measures $\mu$ (see the examples in \cite{werner}). In \cite{kadets} the following 
remarkable result was proved: if $X$ has the DP, then $\norm{\text{id}+T}=1+\norm{T}$ actually 
holds for all weakly compact operators on $X$.\par 
According to \cite{kadets}*{Lemma 2}, $X$ has the DP if and only if for every $x\in S_X$, 
every $x^*\in S_{X^*}$ and all $\eps>0$ there exists $y\in S_X$ such that $x^*(y)\geq 1-\eps$
and $\norm{x+y}\geq 2-\eps$.\par
Thus a test family for the Daugavet property in $X$ is given by
\begin{equation*} 
F_{\eps,U}(\mathbf{x},\mathbf{x}^*,\mathbf{y},\mathbf{y}^*):=
\max\set*{1-x_1^*(y_1),2-\norm{x_1+y_1}}
\end{equation*}
for all $U\in \U(X)$, $\eps>0$, $\mathbf{x}=(x_1,\dots,x_n), \mathbf{y}=(y_1,\dots,y_m)\in B_U^{\text{fin}}$
and all $\mathbf{x}^*=(x_1^*,\dots,x_k^*), \mathbf{y}^*\in B_{U^*}^{\text{fin}}$ (again the conditions (i)--(v)
are easily verified).\par 
The following weaker version of the DP was introduced in \cite{martin}: a real or complex Banach space $X$ is said to have 
the alternative Daugavet property (ADP) if $\max_{\omega\in \T}\norm{\text{id}+\omega T}=1+\norm{T}$ holds for every rank-one 
operator $T$ on $X$, where $\T:=\set*{\omega\in \K:|\omega|=1}$.\par
Again it was proved in \cite{martin} that the above equality holds for all weakly compact opertaors if it holds 
for all rank-one operators. It was also proved in \cite{martin} that $X$ has the ADP if and only if for every $\eps>0$, 
every $x\in S_X$ and every $x^*\in S_{X^*}$ there is some $y\in S_X$ such that $\operatorname{Re}x^*(y)\geq 1-\eps$ and 
$\max_{\omega\in \T}\norm{y+\omega x}\geq 2-\eps$.\par 
We can thus define a test family for the ADP in $X$ as follows:
\begin{equation*} 
F_{\eps,U}(\mathbf{x},\mathbf{x}^*,\mathbf{y},\mathbf{y}^*):=
\max\set*{1-\operatorname{Re}x_1^*(y_1),2-\max_{\omega\in \T}\norm{y_1+\omega x_1}}.
\end{equation*}

\subsection{Lush spaces}\label{sub:lush}
Next we consider the class of lush Banach spaces which was introduced in 
\cite{boyko1} (in connection with the study of the numerical index of Banach spaces).
A Banach space $X$ is called lush provided that for every $\eps>0$ and
all $x_1,x_2\in S_X$ there exists a functional $y^*\in S_{X^*}$ such that
$x_1\in S(y^*,\eps)$ and $\dist{x_2}{\aco{S(y^*,\eps)}}<\eps$, where $\aco$ denotes 
the absolutely convex hull and $d$ is the usual inf-distance.\par
For example, if $K$ is a compact Hausdorff space, then $C(K)$, and more generally every so called
$C$-rich subspace of $C(K)$, is lush (see \cite{boyko1}).\par 
We can define a test family for lushness in $X$ by
\begin{equation*} 
F_{\eps,U}(\mathbf{x},\mathbf{x}^*,\mathbf{y},\mathbf{y}^*):=
\max\set*{1-y_1^*(x_1),\dist{x_2}{\aco{S(y_1^*,\eps)}}}
\end{equation*}
for $\eps>0$, $U\in \U(X)$, $\mathbf{x}=(x_1,\dots,x_n),\mathbf{y}\in B_U^{\text{fin}}$ and 
$\mathbf{x}^*,\mathbf{y}^*=(y_1^*,\dots,y_m^*)\in B_{U^*}^{\text{fin}}$
(where we set $x_2:=x_1$ if $n=1$ and $\dist{x_2}{\aco{S(y_1^*,\eps)}}:=2$ if $\norm{y_1^*}<1$). 
The conditions (i)--(v) in Definition \ref{def:test} are easily verified.\par
In \cite{tan} the following related notion was introduced: the space $X$ 
is called generalised lush (GL) if for every $x\in S_X$ and every $\eps>0$ 
there is some functional $y^*\in S_{X^*}$ such that $x\in S(y^*,\eps)$ and
$\dist{z}{S(y^*,\eps)}+\dist{z}{-S(y^*,\eps)}<2+\eps$ for every $z\in S_X$.\par
It was shown in \cite{tan} that every separable lush space is GL, and that $\R^2$ equipped
with the hexagonal norm $\norm{(a,b)}=\max\set*{\abs{b},\abs{a}+1/2\abs{b}}$ is GL but not lush.
It is not known whether every nonseparable lush space is GL.\par 
The main result in \cite{tan} is that every GL-space $X$ has the Mazur-Ulam property (MUP), i.\,e. if 
$Y$ is any Banach space and $T:S_X \rightarrow S_Y$ is a surjective isometry, then $T$ can be extended 
to an isometric isomorphism between $X$ and $Y$.\par
It is not obvious whether the property GL can be described via test families. However, there is 
the following (at least formally) weaker version of GL-spaces: $X$ is said to have the property $(**)$ 
if for all $x_1,x_2\in S_X$ and each $\eps>0$ one can find $y^*\in S_{X^*}$ such that $x_1\in S(y^*,\eps)$ 
and $\dist{x_2}{S(y^*,\eps)}+\dist{x_2}{-S(y^*,\eps)}<2+\eps$.\par
This notion was introduced in the author's paper \cite{hardtke} (with the help of an anonymous referee)
and the following observations were made:
\begin{enumerate}[(a)]
\item Every lush space has property $(**)$.
\item For separable spaces, $(**)$ is equivalent to GL.
\item Every space with property $(**)$ has the MUP.
\end{enumerate}
A test family for $(**)$ in $X$ can be defined by 
\begin{equation*} 
F_{\eps,U}(\mathbf{x},\mathbf{x}^*,\mathbf{y},\mathbf{y}^*):=
\max\set*{1-y_1^*(x_1),\dist{x_2}{S(y_1^*,\eps)}+\dist{x_2}{-S(y_1^*,\eps)}-2}
\end{equation*}
for $\eps>0$, $U\in \U(X)$, $\mathbf{x}=(x_1,\dots,x_n),\mathbf{y}\in B_U^{\text{fin}}$ and 
$\mathbf{x}^*,\mathbf{y}^*=(y_1^*,\dots,y_m^*)\in B_{U^*}^{\text{fin}}$, where $x_2:=x_1$ if $n=1$
and $\dist{x_2}{S(y_1^*,\eps)}:=\dist{x_2}{-S(y_1^*,\eps)}:=2$ if $\norm{y_1^*}<1$.

\section{Main result}\label{sec:mainresult}
Given a complete, $\sigma$-finite measure space $(S,\A,\mu)$, a K\"othe function space $E$
over $(S,\A,\mu)$ and pairwise disjoint sets $A_1,\dots,A_N\in \A$ with 
$0<\mu(A_i)<\infty$ for $i=1,\dots,N$, we define 
\begin{equation*}
\norm{(a_1,\dots,a_N)}_{E(A_1,\dots,A_N)}:=\norm*{\sum_{i=1}^N\frac{a_i}{\norm{\chi_{A_i}}}_E\chi_{A_i}}_E \ \ \forall (a_1,\dots,a_N)\in \R^N.
\end{equation*} 
Then $\norm{\cdot}_{E(A_1,\dots,A_N)}$ is an absolute, normalised norm on $\R^N$.
If $p\in [1,\infty]$ and $E=L^p(\mu)$, then this norm coincides with the usual $p$-norm on $\R^N$, 
regardless of the choice of $A_1,\dots,A_N$.\par 
For a Banach space $X$, we denote by $E(A_1,\dots,A_N,X)$ the $N$-fold absolute sum
of $X$ with respect to $\norm{\cdot}_{E(A_1,\dots,A_N)}$.\par
The following theorem is the main result of this paper.

\begin{theorem}\label{thm:reduction}
Let $(S,\A,\mu)$ be a complete, $\sigma$-finite measure space and $E$ a K\"othe function space over $(S,\A,\mu)$. 
Suppose that $X$ is a Banach space such that the simple functions are dense in $E(X)$ and 
$E(A_1,\dots,A_N,X)\in \E$ for every $N\in \N$ and all pairwise disjoint sets $A_1,\dots,A_N\in \A$ 
with $0<\mu(A_i)<\infty$ for each $i$. Suppose further that there exists a test family for $\E$ in $E(X)$. 
Then $E(X)\in \E$.	
\end{theorem}

\begin{Proof}
Let $(F_{\eps,U})_{\eps>0,U\in \U(E(X))}$ be a test family for $\E$ in $E(X)$.
Let $\mathbf{f}=(f_1,\dots,f_n)\in S_{E(X)}^n$, $\Phi=(\varphi_1,\dots,\varphi_m)\in S_{E(X)^*}^m$ 
and $\eps>0$. Choose $\delta>0$ such that:
\begin{enumerate}[(a)]
\item For all $\mathbf{y}\in B_{E(X)}^{\text{fin}}$, all $\mathbf{y}^*\in B_{E(X)^*}^{\text{fin}}$ 
and all $\mathbf{z}\in B_{E(X)}^n$ with $\norm{\mathbf{f}-\mathbf{z}}_{\infty}\leq\delta$ we have
\begin{equation*}
\abs{F_{\eps,E(X)}(\mathbf{f},\Phi,\mathbf{y},\mathbf{y}^*)-F_{\eps,E(X)}(\mathbf{z},
\Phi,\mathbf{y},\mathbf{y}^*)}\leq\frac{\eps}{2}. 
\end{equation*}
\item For every $U\in \U(E(X))$, for all $\mathbf{x}\in B_U^n$, all 
$\mathbf{y}\in B_U^{\text{fin}}$, every $\mathbf{y}^*\in B_{U^*}^{\text{fin}}$ and all 
$\mathbf{x}^*,\mathbf{z}^*\in B_{U^*}^m$ with $\norm{\mathbf{x}^*-\mathbf{z}^*}_{\infty}\leq\delta$ 
we have 
\begin{equation*} 
\abs{F_{\eps,U}(\mathbf{x},\mathbf{x}^*,\mathbf{y},\mathbf{y}^*)-
F_{\eps,U}(\mathbf{x},\mathbf{z}^*,\mathbf{y},\mathbf{y}^*)}\leq \frac{\eps}{4c},
\end{equation*}
where $c$ is the constant from Definition \ref{def:test} (iii).
\end{enumerate}
(This is possible because of (iv) and (v) in Definition \ref{def:test}).\par 
Put $\tilde{\eps}:=\min\set*{\eps,\eps/4c}$.\par 
We can find simple functions $h_1,\dots,h_n\in E(X)$ such that $\norm{h_i}_{E(X)}=1$ 
and  $\norm{f_i-h_i}_{E(X)}\leq\delta$ for all $i=1,\dots,n$.\par 
Also, there are simple functions $g_1,\dots,g_m\in E(X)$ with $\norm{g_j}_{E(X)}=1$
and $\abs{\varphi_j(g_j)}\geq 1-\delta$ for all $j=1,\dots,m$.\par
Fix pairwise disjoint sets $A_1,\dots,A_N\in \A$ with $0<\mu(A_i)<\infty$ such that 
each $h_i$ and each $g_j$ belongs to the subspace
\begin{equation*} 
U:=\set*{\sum_{k=1}^Nx_k\chi_{A_k}:x_1,\dots,x_N\in X}\ssq E(X).
\end{equation*}
By considering the map $T:E(A_1,\dots,A_N,X) \rightarrow U$ defined by
\begin{equation*}
T(x_1,\dots,x_N):=\sum_{k=1}^N\frac{x_k}{\norm{\chi_{A_k}}_E}\chi_{A_k}, 
\end{equation*}
we see that $U$ is isometrically isomorphic to $E(A_1,\dots,A_N,X)$.\par 
By assumption we have $E(A_1,\dots,A_N,X)\in \E$, thus $U\in \E$.\par
Since $g_j\in S_U$ we have $1\geq\norm{\varphi_j|_U}\geq 1-\delta$ for each $j$.
Hence $\psi_j:=\varphi_j|_U/\norm{\varphi_j|_U}\in S_{U^*}$ with 
\begin{equation}\label{eq:psi}
\norm{\psi_j-\varphi_j|_U}=\abs{1-\norm{\varphi_j|_U}}\leq\delta \ \ \forall j=1,\dots,m.
\end{equation}
Put $\Psi:=(\psi_1,\dots,\psi_m)\in S_{U^*}^m$ and $\mathbf{h}=(h_1,\dots,h_n)\in S_U^n$. Since $U\in \E$ 
we can find $\mathbf{u}=(u_1,\dots,u_l)\in S_U^{\text{fin}}$ and $\mathbf{u}^*=(u_1^*,\dots,u_s^*)\in S_{U^*}^{\text{fin}}$ such that 
$F_{\tilde{\eps},U}(\mathbf{h},\Psi,\mathbf{u},\mathbf{u}^*)\leq \tilde{\eps}$.\par 
Because of $\tilde{\eps}\leq \eps$ and (ii) in Definition \ref{def:test}, it follows that 
$F_{\eps,U}(\mathbf{h},\Psi,\mathbf{u},\mathbf{u}^*)\leq \tilde{\eps}$.\par
Then (b) and \eqref{eq:psi} imply 
$F_{\eps,U}(\mathbf{h},\Phi|_U,\mathbf{u},\mathbf{u}^*)\leq \tilde{\eps}+\eps/4c\leq \eps/2c$.\par 
By the Hahn-Banach theorem there are functionals $\omega_1,\dots,\omega_s\in S_{E(X)^*}$ 
such that $\omega_i|_U=u_i^*$ for $i=1,\dots,s$. Let $\Omega:=(\omega_1,\dots,\omega_s)$.\par
Now it follows from (iii) in Definition \ref{def:test} that 
$F_{\eps,E(X)}(\mathbf{h},\Phi,\mathbf{u},\Omega)\leq\eps/2$.\par 
Since $\norm{\mathbf{f}-\mathbf{h}}_{\infty}\leq\delta$, (a) imlpies 
$F_{\eps,E(X)}(\mathbf{f},\Phi,\mathbf{u},\Omega)\leq\eps/2+\eps/2=\eps$
and the proof is finished.
\end{Proof}

Every K\"othe function space $E$ is a Banach lattice in its natural 
ordering ($f\leq g$ if and only if $f(s)\leq g(s)$ for a.\,e. $s\in S$).
It is well-known that if $(E,\leq)$ is order continuous, then for every 
Banach space $X$ the simple functions lie dense in $E(X)$. This includes 
in particular the case of $L^p$-spaces for $1\leq p<\infty$. So from the 
above theorem we obtain the following corollary ($\ell^p_N(X)$ denotes the
$N$-fold $p$-sum of $X$).

\begin{corollary}\label{cor:reductionp}
Let $(S,\A,\mu)$ be a complete, $\sigma$-finite measure space and $1\leq p<\infty$. If $X$ is 
a Banach space such that $\ell^p_N(X)\in \E$ for every $N\in \N$ and there exists a 
test family for $\E$ in $L^p(\mu,X)$, then $L^p(\mu,X)\in \E$.
\end{corollary}

In the case $p=\infty$, it is well-known that one still has the density of 
$\set*{f\in L^{\infty}(\mu,X):\text{ran}(f)\ \text{is\ countable}}$ in $L^{\infty}(\mu,X)$,
where $\text{ran}(f)$ denotes the range of $f$. Thus one can prove the following Theorem 
in an analogous way to the proof of Theorem \ref{thm:reduction} (we omit the details).
\begin{theorem}\label{thm:Linfty}
Let $(S,\A,\mu)$ be a complete, $\sigma$-finite measure space. If $X$ is a 
Banach space such that $\ell_N^{\infty}(X)\in \E$ for every $N\in \N$ and 
$\ell^{\infty}(X)\in\E$ and there exists a test family for $\E$ in $L^{\infty}(\mu,X)$, 
then $L^{\infty}(\mu,X)\in \E$. 	
\end{theorem}
Here $\ell^{\infty}(X)$ stands for $\big[\bigoplus_{n\in \N}X\big]_{\ell^{\infty}}$.\par 
We also have a reduction result for the case of infinite absolute sums to finite sums,
which reads as follows.
\begin{proposition}\label{prop:sums}
Let $I$ be an index set and $E$ a subspace of $\R^I$ endowed with an absolute,
normalised norm such that $\operatorname{span}\set*{e_i:i\in I}$ is dense in $E$. 
Let $(X_i)_{i\in I}$ be a family of Banach spaces such that $\big[\bigoplus_{i\in J}X_i\big]_E\in \E$ 
for every nonempty, finite subset $J\ssq I$. If there is a test family for $\E$ in 
$\big[\bigoplus_{i\in I}X_i\big]_E$, then $\big[\bigoplus_{i\in I}X_i\big]_E\in \E$.
\end{proposition}
The notation $\big[\bigoplus_{i\in J}X_i\big]_E$ means that all summands with index in $I\sm J$ are $\set*{0}$.
The proof is similar to the one of Theorem \ref{thm:reduction} and will therefore be omitted.\par
As an immediate consequence of Proposition \ref{prop:sums} we get the following results for $p$-sums
and $c_0$-sums.
\begin{corollary}\label{cor:sum-p}
If $I$ is any index set, $1\leq p<\infty$, $(X_i)_{i\in I}$ is a family of Banach spaces 
such that $\big[\bigoplus_{i\in J}X_i\big]_p\in \E$ for every nonempty, finite subset $J\ssq I$, 
and there exists a test family for $\E$ in $\big[\bigoplus_{i\in I}X_i\big]_p$, then 
$\big[\bigoplus_{i\in I}X_i\big]_p\in \E$.
\end{corollary}

\begin{corollary}\label{cor:c0-sum}
If $I$ is any index set, $(X_i)_{i\in I}$ is a family of Banach spaces such that
$\big[\bigoplus_{i\in J}X_i\big]_{\infty}\in \E$ for every nonempty, finite subset 
$J\ssq I$, and there exists a test family for $\E$ in $\big[\bigoplus_{i\in I}X_i\big]_{c_0}$, 
then $\big[\bigoplus_{i\in I}X_i\big]_{c_0}\in \E$.
\end{corollary}

\section{Applications}\label{sec:appl}
In this section we will apply the abstract results to the examples discussed earlier. This
will yield some new results as well as some alternative proofs of already known results.\par 
We first collect what is known about sums of octahedral spaces and their relatives. The following 
results were proved in \cite{haller3}: if $X$ and $Y$ are real Banach spaces, then
\begin{enumerate}[(a)]
\item $X$ or $Y$ is LOH/WOH/OH $\Rightarrow$ $X\oplus_1 Y$ is LOH/WOH/OH,
\item $X$ and $Y$ are LOH/WOH $\Rightarrow$ $X\oplus_p Y$ is LOH/WOH for every $p\in (1,\infty]$,
\item $X$ and $Y$ are OH $\Rightarrow$ $X\oplus_{\infty} Y$ is OH,
\item For $p\in (1,\infty)$ $X\oplus_p Y$ is never OH.
\end{enumerate}

In \cite{abrahamsen2} the following generalisation was obtained: if $I$ is any index set and $E$ a 
subspace of $\R^I$ with an absolute, normalised norm, and $(X_i)_{i\in I}$ is a family of LOH spaces,
then $\big[\bigoplus_{i\in I}X_i\big]_E$ is also LOH. If each $X_i$ is WOH and moreover 
$\text{span}\set*{e_i:i\in I}$ is dense in $E$, then $\big[\bigoplus_{i\in I}X_i\big]_E$ is also WOH.\par
It is also easily checked that $\ell^{\infty}(X)$ is OH whenever $X$ is OH (the proof is analogous to the proof
of (c) above that was given in \cite{haller3}).\par 
Combining all this with our Theorems \ref{thm:reduction} resp. \ref{thm:Linfty} and the fact that OH, WOH 
and LOH can be described by test families (see Section \ref{sec:examples}), we obtain the following results.
\begin{theorem}\label{thm:LOH-WOH}
If $(S,\A,\mu)$ is a complete, $\sigma$-finite measure space, $E$ a K\"othe function space over $(S,\A,\mu)$ 
and $X$ an LOH/WOH space such that the simple functions are dense in $E(X)$ (for instance, if $E$ is order continuous), 
then $E(X)$ is also LOH/WOH.\par 
In particular, if $p\in [1,\infty)$ and $X$ is LOH/WOH, then so is $L^p(\mu,X)$.\par 
Also, $L^{\infty}(\mu,X)$ is LOH if $X$ is LOH.
\end{theorem}

\begin{proposition}\label{prop:OH}
If $(S,\A,\mu)$ is a complete, $\sigma$-finite measure space and $X$ is an OH space, then
$L^1(\mu,X)$ and $L^{\infty}(\mu,X)$ are also OH.
\end{proposition}
This result is not optimal. In fact, it is not difficult to see that $L^1(\mu,X)$ is OH for 
{\it any} Banach space $X$ (provided that $L^1(\mu)$ is infinite-dimensional), see the examples
at the end of \cite{langemets1}.\par 
Now we turn to the diameter-two-properties. In \cite{haller3} the following results were derived 
via duality from the corresponding results on octahedrality in direct sums.
\begin{enumerate}[(a)]
\item $X$ or $Y$ has the LD2P/D2P/SD2P $\Rightarrow$ $X\oplus_{\infty} Y$ has the LD2P/D2P/SD2P,
\item $X$ and $Y$ have the LD2P/D2P $\Rightarrow$ $X\oplus_p Y$ has the LD2P/D2P for every $p\in [1,\infty)$,
\item $X$ and $Y$ have the SD2P $\Rightarrow$ $X\oplus_1 Y$ has the SD2P,
\item For $p\in (1,\infty)$ $X\oplus_p Y$ never has the SD2P.
\end{enumerate}
All these results have been known before (they are scattered in \cites{abrahamsen, acosta, becerra-guerrero, haller, lopez-perez},
see \cite{haller3} for a detailed account), but the previous proofs were based on different methods. In \cite{acosta} it was shown
that the LD2P and the D2P are stable under sums with respect to an arbitrary absolute norm.\par
Since LD2P, D2P and SD2P can be described by test families (see Section \ref{sec:examples}), we obtain the following stability
result from Theorem \ref{thm:reduction}.
\begin{theorem}\label{thm:diam2}
Let $(S,\A,\mu)$ be a complete, $\sigma$-finite measure space, $E$ a K\"othe function space over $(S,\A,\mu)$ and $X$ a Banach 
space such that the simple functions are dense in $E(X)$ (for instance, if $E$ is order continuous). If $X$ has the LD2P/D2P, then
$E(X)$ also has the LD2P/D2P.\par 
In particular, if $p\in [1,\infty)$ and $X$ has the LD2P/D2P, then so does $L^p(\mu,X)$.\par
Further, if $X$ has the SD2P, then $L^1(\mu,X)$ also has the SD2P.
\end{theorem}

In \cite{acosta} it was already proved that $L^p(\mu,X)$ has the D2P whenever $1\leq p <\infty$, $\mu$ is a finite measure and
$X$ has the D2P (this proof also uses simple functions). Also, for the special case $p=1$, better results are already known, for 
instance, it has been proved in \cite{becerra-guerrero}*{Theorem 2.13} that for a finite measure $\mu$ the space $L^1(\mu,X)$ has 
the D2P if and only if $X$ has the D2P or $\mu$ has no atoms (and $L^{\infty}(\mu,X)$ has the D2P if and only if $L^{\infty}(\mu)$ 
is infinite-dimensional or $X$ has the D2P).\par
Even more, it is known that the Daugavet property implies the SD2P (see \cite{abrahamsen}*{Theorem 4.4}) and that
$L^1(\mu,X)$ and $L^{\infty}(\mu,X)$ have the Daugavet property for any atomless measure $\mu$ and any Banach space $X$
(\cite{werner}, see the discussion for the DP below).\par
Also, if $X$ or $Y$ has the LD2P, then so does $X\hat{\otimes}_{\pi}Y$ (see \cite{abrahamsen}*{Theorem 2.7}) and if $X$ 
and $Y$ have the SD2P, then so does $X\hat{\otimes}_{\pi}Y$ (see \cite{becerra-guerrero5}), where $\hat{\otimes}_{\pi}$ denotes the projective 
tensor product, and it is well-known that $L^1(\mu,X)=L^1(\mu)\hat{\otimes}_{\pi}X$. For more information on octahedrality and related 
properties in tensor products see also \cites{langemets1,langemets2}.\par 
For AOH spaces, the following equivalent characterisation can be proved:
$X$ is AOH if and only if for every $n\in \N$, all $x_1,\dots,x_n\in S_X$ and each $\eps>0$ there is some $y\in S_X$
such that 
\begin{equation*}
\max\set*{\norm{x_i+ty},\norm{x_i-ty}}\geq (1-\eps)(1+t) \ \ \forall t>0, \forall i\in \set*{1,\dots,n}.
\end{equation*}
The proof is analogous to the proof of the corresponding characterisation for octahedral spaces in \cite{haller3}
and will therefore be skipped.\par
Using this characterisation, one can show that $X\oplus_1 Y$ is AOH if $X$ or $Y$ is AOH and that $X\oplus_{\infty} Y$
is AOH if $X$ and $Y$ are AOH. The latter result also extends to $\ell^{\infty}(X)$. Again the proofs are analogous 
to the ones for the corresponding results on OH spaces in \cite{haller3} and thus we will skip them.\par
Using our Theorems \ref{thm:reduction} and \ref{thm:Linfty} we can now obtain the following result.
\begin{proposition}\label{prop:AOH}
If $(S,\A,\mu)$ is a complete, $\sigma$-finite measure space and $X$ is an AOH space, then
$L^1(\mu,X)$ and $L^{\infty}(\mu,X)$ are also AOH.
\end{proposition}
Again, if $L^1(\mu)$ is infinite-dimensional, then $L^1(\mu,X)$ is even OH for {\it any} Banach space $X$ (\cite{langemets1}).\par 
Concerning sums of ASQ and LASQ spaces, the following was proved in \cite{abrahamsen2}: if $I$ is any index 
set and $E$ a subspace of $\R^I$ with an absolute, normalised norm, and $(X_i)_{i\in I}$ is a family of LASQ
spaces, then $\big[\bigoplus_{i\in I}X_i\big]_E$ is also LASQ. Further, $X\oplus_{\infty} Y$ is ASQ/LASQ if 
and only if $X$ or $Y$ is ASQ/LASQ. Analogously to the proof of the ``if'' part in \cite{abrahamsen2} one can 
show that $\ell^{\infty}(X)$ is ASQ/LASQ whenever $X$ is ASQ/LASQ.\footnote{It has also been proved in 
\cite{abrahamsen2} that for $p\in [1,\infty)$ the sum $X\oplus_p Y$ is never ASQ.}\par
If we combine these facts with Theorem \ref{thm:reduction} resp. \ref{thm:Linfty} and the fact that ASQ and LASQ
can be expressed in terms of test families (Section \ref{sec:examples}), we obtain the following stability result.
\begin{theorem}\label{ASQ-LASQ}
If $(S,\A,\mu)$ is a complete, $\sigma$-finite measure space, $E$ a K\"othe function space over $(S,\A,\mu)$ 
and $X$ an LASQ space such that the simple functions are dense in $E(X)$ (for instance, if $E$ is order continuous), 
then $E(X)$ is also LASQ.\par
In particular, if $p\in [1,\infty)$ and $X$ is LASQ, then so is $L^p(\mu,X)$.\par 
Moreover, $L^{\infty}(\mu,X)$ is ASQ/LASQ whenever $X$ is ASQ/LASQ.
\end{theorem}

Now we consider spaces with the Daugavet property. It has been shown in \cite{khalil} that
$L^1([0,1],X)$ and $L^{\infty}([0,1],X)$ have the DP if $X$ has it.
More generally, $L^1(\mu,X)$ has the DP for every {\it atomless} measure $\mu$ and 
{\it every} Banach space $X$, see \cite{werner}*{p.81}.\par
In \cite{wojtaszczyk} it was already proved that the $\ell^1$- and $\ell^{\infty}$-sum of any (finite or infinite) 
sequence of Banach spaces with the Daugavet property has again the Daugavet property (the Daugavet property for 
weakly compact operators was considered in \cite{wojtaszczyk}, but this is equivalent to considering just rank-one 
operators by \cite{kadets}*{Theorem 2.3}). In \cite{kadets} a different proof for the stability of the DP by finite 
or infinite $\ell^1$- and $c_0$-sums has been given.\footnote{The cases of infinite sums are reduced to the corresponding 
finite sums by a density argument, similar to the general reduction results for sums that we have stated in Section \ref{sec:mainresult}.}\par
Putting everything together, the following characterisation was obtained in \cite{martin00}*{Remark 9}: $L^1(\mu,X)$ has 
the DP if and only if $X$ has the DP or $\mu$ has no atoms. Likewise, $L^{\infty}(\mu,X)$ has the DP if and only if $X$ has 
the DP or $\mu$ has no atoms (see \cite{martin0}).\par
Analogous results also hold for the alternative Daugavet property: the space $L^1(\mu,X)$ has the ADP if and only if $X$ has the ADP
or $\mu$ has no atoms if and only if $L^{\infty}(\mu,X)$ has the ADP (see \cite{martin}). Also, the ADP is stable under arbitrary 
$\ell^1$-, $c_0$- and $\ell^{\infty}$-sums (see again \cite{martin}).\par 
Using the stability results for sums and our Theorems \ref{thm:reduction} and \ref{thm:Linfty}, we obtain an alternative 
proof of the following known result.\par
\begin{theorem}\label{thm:daugavet}
If $(S,\A,\mu)$ is a complete, $\sigma$-finite measure space and $X$ a Banach space with the DP/ADP,
then $L^1(\mu,X)$ and $L^{\infty}(\mu,X)$ also have the DP/ADP.
\end{theorem}

Concerning lush spaces, the following has been proved in \cite{boyko2}: if $\norm{\cdot}_E$ is an absolute norm on $\R^n$, 
then the sum of every collection $X_1,\dots,X_n$ of lush spaces with respect to $\norm{\cdot}_E$ is again lush if and only if 
$(\R^n,\norm{\cdot}_E)$ is lush.\par 
It was also proved in \cite{boyko2} that the $\ell^1$-, $c_0$- and $\ell^{\infty}$-sums of any family $(X_i)_{i\in I}$ of 
lush spaces are again lush.\footnote{Also here the cases of $\ell^1$- and $c_0$-sums are reduced to the corresponding 
finite sums (cf. footnote 3).}\par
Very recently, the following stability result has been proved in \cite{kadets2}*{Corollaries 8.9 and 8.12}.
\begin{theorem}[\cite{kadets2}]\label{thm:lush}
Let $(S,\A,\mu)$ be a $\sigma$-finite measure space and $X$ a Banach space. Then $L^{\infty}(\mu,X)$ is lush if and only if $X$
is lush if and only if $L^1(\mu,X)$ is lush.
\end{theorem}
In fact, even more general results are proved in \cite{kadets2} for so called lush operators.\par 
If we use instead the above-mentioned results on sums of lush spaces in combination with our Theorems \ref{thm:reduction} 
and \ref{thm:Linfty}, we obtain an alternative proof for the fact that lushness of $X$ is sufficient for lushness of 
$L^1(\mu,X)$ and $L^{\infty}(\mu,X)$ (the proofs in \cite{kadets2} did not use a reduction to sums, but they also used 
the density of the simple functions (resp. functions with countable range) in $L^1(\mu,X)$ (resp. $L^{\infty}(\mu,X)$).\par 
Let us now turn to generalised lushness. It was proved in  \cite{tan} that the property GL is stable under arbitrary 
$\ell^1$-, $c_0$- and $\ell^{\infty}$-sums. The same results also hold for the property $(**)$, with completely analogous proofs.\par
Now we can apply our Theorems \ref{thm:reduction} and \ref{thm:Linfty} to obtain the following result.\par 
\begin{theorem}\label{thm:starstar}
If $(S,\A,\mu)$ is a complete, $\sigma$-finite measure space and $X$ a Banach space with property $(**)$, then
$L^1(\mu,X)$ and $L^{\infty}(\mu,X)$ also have the property $(**)$.
\end{theorem}
We recall (see Subsection \ref{sub:lush}) that $(**)$ implies the MUP and $(**)$ is equivalent to GL for separable spaces, but it is 
not known whether this equivalence is true in general nor if it is in general possible to describe the property GL by test families. 
Thus we cannot apply our general reduction theorems directly to GL-spaces. However, it is still possible to show that GL is stable 
with respect to $L^1$-Bochner spaces by a similar proof technique. This is carried out in the next section.

\section{GL-spaces}\label{sec:GL}
Here we show directly that $L^1(\mu,X)$ is GL whenever $X$ is GL. The argument is similar to the proof for $\ell^1$-sums
given in \cite{tan}, in combination with an approximation by simple functions.

\begin{theorem}\label{thm:GL}
Let $(S,\A,\mu)$ be a complete, $\sigma$-finite measure space. If $X$ is a GL-space,
then so is $L^1(\mu,X)$.
\end{theorem}

\begin{Proof} 
Let $f\in L^1(\mu,X)$ with $\norm{f}_1=1$ and let $\eps\in (0,1)$. Choose $\eta\in (0,1)$ such that
\begin{align*}
&(2+\eps/2)(1+\eta)+4\eta<2+\eps, \\
&(1-\eps/2)(1-\eta)-\eta>1-\eps, \\
&(1-\eps/2)\frac{1-\eta}{1+\eta}>1-\eps. 
\end{align*} 
We can find a simple function $g$ on $S$ such $\norm{f-g}_1\leq \eta$.
Write $g=\sum_{i=1}^Nx_i\chi_{A_i}$ with pairwise disjoint sets $A_1,\dots,A_N\in \A$
and $x_1,\dots,x_N\in X$.\par 
Since $X$ is GL, we can find functionals $x_1^*,\dots,x_N^*\in S_{X^*}$ such that
$x_i^*(x_i)\geq (1-\eps/2)\norm{x_i}$ and
\begin{equation}\label{eq:GL} 
\dist{y}{S(x_i^*,\eps/2)}+\dist{y}{-S(x_i^*,\eps/2)}<2+\frac{\eps}{2} \ \ \forall y\in S_X.
\end{equation}
Let $h=\sum_{i=1}^Nx_i^*\chi_{A_i}$ and 
$\varphi(v)=\int_Sh(s)(v(s))\,\text{d}\mu(s)$ for $v\in L^1(\mu,X)$. Then 
$\varphi\in L^1(\mu,X)^*$ with $\norm{\varphi}=1$.\par 
We further have
\begin{equation*} 
\varphi(g)=\sum_{i=1}^N\int_{A_i}x_i^*(x_i)\,\text{d}\mu(s)\geq 
(1-\eps/2)\sum_{i=1}^N\int_{A_i}\norm{x_i}\,\text{d}\mu(s)=(1-\eps/2)\norm{g}_1.
\end{equation*} 
Since $\norm{f-g}_1\leq \eta$ and $\norm{f}_1=1$, it follows that 
$\varphi(f)\geq \varphi(g)-\eta\geq (1-\eps/2)(1-\eta)-\eta$. Thus the choice of $\eta$ 
implies $f\in S(\varphi,\eps)$.\par
Now take any function $w\in L^1(\mu,X)$ with $\norm{w}_1=1$. There exists a simple
function $\tilde{w}$ on $S$ such that $\norm{w-\tilde{w}}_1\leq \eta$. Write $\tilde{w}=\sum_{j=1}^My_j\chi_{B_j}$ 
with pairwise disjoint sets $B_1,\dots,B_M\in \A$ and $y_1,\dots,y_M\in X$.\par 
We put $C_{ij}:=A_i\cap B_j$ for $(i,j)\in I:=\set*{1,\dots,N}\times \set*{1,\dots,M}$. 
By \eqref{eq:GL} we can find, for each pair $(i,j)\in I$, vectors $u_{ij},v_{ij}\in B_X$ such that
$x_i^*(u_{ij})>1-\eps/2$, $-x_i^*(v_{ij})>1-\eps/2$ and 
\begin{equation}\label{eq:GL2}
\norm{y_j-\norm{y_j}u_{ij}}+\norm{y_j-\norm{y_j}v_{ij}}\leq(2+\eps/2)\norm{y_j}.
\end{equation}
Let $u=\sum_{(i,j)\in I}\norm{y_j}u_{ij}\chi_{C_{ij}}$ and    
$v=\sum_{(i,j)\in I}\norm{y_j}v_{ij}\chi_{C_{ij}}$.\par 
Since $\norm{u_{ij}}\leq 1$ we have $\norm{u(s)}\leq \norm{y_j}=\norm{\tilde{w}(s)}$
for all $s\in C_{ij}$ and all $(i,j)\in I$. Hence $\norm{u}_1\leq \norm{\tilde{w}}_1
\leq \norm{w-\tilde{w}}_1+\norm{w}_1\leq 1+\eta$. Analogously, one can see that 
$\norm{v}_1\leq 1+\eta$.\par 
Thus we have $\tilde{u}:=u/(1+\eta)\in B_{L^1(\mu,X)}$ and 
$\tilde{v}:=v/(1+\eta)\in B_{L^1(\mu,X)}$.\par 
We further have
\begin{align*} 
&\varphi(\tilde{u})=\frac{1}{1+\eta}\sum_{(i,j)\in I}\int_{C_{ij}}x_i^*(u_{ij})\norm{y_j}\,\text{d}\mu(s)\geq 
\frac{1-\eps/2}{1+\eta}\sum_{(i,j)\in I}\int_{C_{ij}}\norm{y_j}\,\text{d}\mu(s) \\
&=\norm{\tilde{w}}_1\frac{1-\eps/2}{1+\eta}\geq (1-\eta)\frac{1-\eps/2}{1+\eta}
>1-\eps.
\end{align*} 
Thus $\tilde{u}\in S(\varphi,\eps)$ and analogously one can show that 
$\tilde{v}\in -S(\varphi,\eps)$.\par
It further follows from \eqref{eq:GL2} that 
\begin{equation*} 
\norm{\tilde{w}(s)-u(s)}+\norm{\tilde{w}(s)-v(s)}\leq (2+\eps/2)\norm{\tilde{w}(s)}
\ \ \forall s\in S.
\end{equation*}  
Hence $\norm{\tilde{w}-u}_1+\norm{\tilde{w}-v}_1\leq (2+\eps/2)\norm{\tilde{w}}_1\leq (2+\eps/2)(1+\eta)$.\par 
Since $\norm{w-\tilde{w}}_1\leq \eta$ we get 
$\norm{w-u}_1+\norm{w-v}_1\leq (2+\eps/2)(1+\eta)+2\eta$.\par 
We also have $\norm{u-\tilde{u}}_1\leq \eta$ and $\norm{v-\tilde{v}}_1\leq \eta$.
Thus $\norm{w-\tilde{u}}_1+\norm{w-\tilde{v}}_1\leq (2+\eps/2)(1+\eta)+4\eta
<2+\eps$ and we are done.
\end{Proof}

\address
\email


\begin{bibdiv}
\begin{biblist}

\bib{abrahamsen}{article}{
  title={Remarks on diameter 2 properties},
  author={Abrahamsen, T.},
  author={Lima, V.},
  author={Nygaard, O.},
  journal={J. Conv. Anal.},
  volume={20},
  date={2013},
  pages={439--452}
  }
  
\bib{abrahamsen2}{article}{
  title={Almost square Banach spaces},
  author={Abrahamsen, T. A.},
  author={Langemets, J.},
  author={Lima, V.},
  journal={J. Math. Anal. Appl.},
  volume={434},
  number={2},
  date={2016},
  pages={1549--1565}
  }
  
\bib{acosta}{article}{
  title={Stability results of diameter two properties},
  author={Acosta, M. D.},
  author={Becerra Guerrero, J.},
  author={L\'opez-P\'erez, G.},
  journal={J. Conv. Anal.},
  volume={22},
  number={1},
  date={2015},
  pages={1--17}
  }
  
\bib{becerra-guerrero}{article}{
  title={Relatively weakly open subsets of the unit ball in function spaces},
  author={Becerra Guerrero, J.},
  author={L\'opez-P\'erez, G.},
  journal={J. Math. Anal. Appl.},
  volume={315},
  date={2006},
  pages={544--554}
  }
  
\bib{becerra-guerrero2}{article}{
  title={Octahedral norms and convex combination of slices in Banach spaces},
  author={Becerra Guerrero, J.},
  author={L\'opez-P\'erez, G.},
  author={Rueda Zoca, A.},
  journal={J. Funct. Anal.},
  volume={266},
  number={4},
  date={2014},
  pages={2424--2435}
  }
    
\bib{becerra-guerrero4}{article}{
  title={Big slices versus big relatively weakly open subsets in Banach spaces},
  author={Becerra Guerrero, J.},
  author={L\'opez-P\'erez, G.},
  author={Rueda Zoca, A.},
  journal={J. Math. Anal. Appl.},
  volume={428},
  date={2015},
  pages={855--865}
  }

\bib{becerra-guerrero5}{article}{
  title={Octahedral norms in spaces of operators},
  author={Becerra Guerrero, J.},
  author={L\'opez-P\'erez, G.},
  author={Rueda Zoca, A.},
  journal={J. Math. Anal. Appl.},
  volume={427},
  date={2015},
  pages={171--184}
  }
  
\bib{becerra-guerrero3}{article}{
  title={Some results on almost square Banach spaces},
  author={Becerra Guerrero, J.},
  author={L\'opez-P\'erez, G.},
  author={Rueda Zoca, A.},
  journal={J. Math. Anal. Appl.},
  volume={438},
  number={2},
  date={2016},
  pages={1030--1040}
  }

\bib{boyko1}{article}{
  title={Numerical index of Banach spaces and duality},
  author={Boyko, K.},
  author={Kadets, V.},
  author={Mart\'{\i}n, M.},
  author={Werner, D.},
  journal={Math. Proc. Cambridge Philos. Soc.},
  volume={142},
  number={1},
  date={2007},
  pages={93--102}
  }
  
\bib{boyko2}{article}{
  title={Properties of lush spaces and applications to Banach spaces with numerical index one},
  author={Boyko, K.},
  author={Kadets, V.},
  author={Mart\'{\i}n, M.},
  author={Mer\'{\i}, J.},
  journal={Studia Math.},
  volume={190},
  date={2009},
  pages={117--133}
  }
  
\bib{deville}{book}{
  title={Smoothness and renormings in Banach spaces},
  author={Deville, R.},
  author={Godefroy, G.},
  author={Zizler, V.},
  publisher={Longman Scientific \& Technical},
  address={Harlow},
  series={Pitman Monographs and Surveys in Pure and Applied Mathematics},
  volume={64},
  date={1993}
  }
  
\bib{godefroy}{article}{
  title={Metric characterization of first Baire class linear forms and octahedral norms},
  author={Godefroy, G.},
  journal={Studia Math.},
  volume={95},
  number={1},
  date={1989},
  pages={1--15}
  }
  
\bib{haller}{article}{
  title={Two remarks on diameter 2 properties},
  author={Haller, R.},
  author={Langemets, J.},
  journal={Proc. Estonian Acad. Sci.},
  volume={63},
  number={1},
  date={2014},
  pages={2--7}
  }
  
\bib{haller2}{article}{
  title={Geometry of Banach spaces with an octahedral norm},
  author={Haller, R.},
  author={Langemets, J.},
  journal={Acta Comment. Univ. Tartuensis Math.},
  volume={18},
  number={1},
  date={2014},
  pages={125-133}
  }
  
\bib{haller3}{article}{
  title={On duality of diameter 2 properties},
  author={Haller, R.},
  author={Langemets, J.},
  author={P\~{o}ldvere, M.},
  journal={J. Conv. Anal.},
  volume={22},
  number={2},
  date={2015},
  pages={465--483}
  }
  
\bib{haller4}{article}{
  title={Rough norms in spaces of operators},
  author={Haller, R.},
  author={Langemets, J.},
  author={P\~{o}ldvere, M.},
  journal={Math. Nachr.},
  date={2017},
  pages={11p.},
  note={doi:10.1002/mana.201600409}
  }
  
\bib{hardtke}{article}{
  title={Some remarks on generalised lush spaces},
  author={Hardtke, J.-D.},
  journal={Studia Math.},
  volume={231},
  number={1},
  date={2015},
  pages={29--44}
  }

\bib{kadets}{article}{
  title={Banach spaces with the Daugavet property},
  author={Kadets, V. M.},
  author={Shvidkoy, R. V.},
  author={Sirotkin, G. G.},
  author={Werner, D.},
  journal={Trans. Amer. Math. Soc.},
  volume={352},
  number={2},
  date={2000},
  pages={855--873}
  }
  
\bib{kadets2}{article}{
  title={Spear operators between Banach spaces},
  author={Kadets, V.},
  author={Mart\'{\i}n, M.},
  author={Mer\'{i}, J.},
  author={P\'erez, A.},
  date={2017},
  pages={114p.},
  note={Preprint, available at \url{http://arxiv.org/abs/1701.02977}}
  }
  
\bib{khalil}{article}{
  title={The Daugavet equation in vector-valued function spaces},
  author={Khalil, R.},
  journal={Panam. Math. J.},
  volume={6},
  number={3},
  date={1996},
  pages={51--53}
  }
  
\bib{kubiak}{article}{
  title={Some geometric properties of the Ces\`aro function spaces},
  author={Kubiak, D.},
  journal={J. Convex Anal.},
  volume={21},
  number={1},
  date={2014},
  pages={189--200}
}
  
\bib{langemets1}{article}{
  title={Almost square and octahedral norms in tensor products of Banach spaces},
  author={Langemets, J.},
  author={Lima, V.},
  author={Rueda Zoca, A.},
  date={2016},
  pages={16p.},
  note={Preprint, available at \url{http://arxiv.org/abs/1602.07090}}
  }
  
\bib{langemets2}{article}{
  title={Octahedral norms in tensor products of Banach spaces},
  author={Langemets, J.},
  author={Lima, V.},
  author={Rueda Zoca, A.},
  date={2016},
  pages={15p.},
  note={Preprint, available at \url{http://arxiv.org/abs/1609.02062}}
  }
  
\bib{lee}{article}{
  title={Polynomial numerical indices of Banach spaces with absolute norms},
  author={Lee, H. J.},
  author={Mart\'{\i}n, M.},
  author={Mer\'{\i}, J.},
  journal={Linear Algebra and its Applications},
  volume={435},
  date={2011},
  pages={400--408}
  }
  
\bib{lin}{book}{
  title={K\"othe-Bochner function spaces},
  author={Lin, P. K.},
  publisher={Birkh\"auser},
  address={Boston-Basel-Berlin},
  date={2004}
  }
  
\bib{lopez-perez}{article}{
  title={The big slice phenomena in $M$-embedded and $L$-embedded spaces},
  author={L\'opez-P\'erez, G.},
  journal={Proc. Amer. Math. Soc.},
  volume={134},
  date={2005},
  pages={273--282}
  }
  
\bib{martin00}{article}{
  title={Numerical index of vector-valued function spaces},
  author={Mart\'{\i}n, M.},
  author={Pay\'a, R.},
  journal={Studia Math.},
  volume={142},
  date={2000},
  pages={269--280}
  }
  

\bib{martin0}{article}{
  title={Numerical index and Daugavet property for $L^{\infty}(\mu,X)$},
  author={Mart\'{\i}n, M.},
  author={Villena, A. R.},
  journal={Proc. Edinburgh Math. Soc.},
  volume={46},
  date={2003},
  pages={415--420}
  }
  
\bib{martin}{article}{
  title={An alternative Daugavet property},
  author={Mart\'{\i}n, M.},
  author={Oikhberg, T.},
  journal={J. Math. Anal. Appl.},
  volume={294},
  number={1},
  date={2004},
  pages={158--180}
  }
  
\bib{tan}{article}{
  title={Generalized-lush spaces and the Mazur-Ulam property},
  author={Tan, D.},
  author={Huang, X.},
  author={Liu, R.},
  journal={Studia Math.},
  volume={219},
  number={2},
  date={2013},
  pages={139--153}
  }
  
\bib{werner}{article}{
  title={Recent progress on the Daugavet property},
  author={Werner, D.},
  journal={Irish Math. Soc. Bulletin},
  volume={46},
  date={2001},
  pages={77--97}
  }
  
\bib{wojtaszczyk}{article}{
  title={Some remarks on the Daugavet equation},
  author={Wojtaszczyk, P.},
  journal={Proc. Amer. Math. Soc.},
  volume={115},
  number={4},
  date={1992},
  pages={1047--1052}
  }
  
\end{biblist}
\end{bibdiv}
\end{document}